\documentclass[11pt]{article}
\usepackage{amsfonts}
\usepackage{amsfonts}
\usepackage{amsfonts}
\usepackage{mathrsfs}
\usepackage{bm}
\usepackage{amssymb,amsmath} 
 \allowdisplaybreaks

\catcode`!=11
\let\!int\int \def\int{\displaystyle\!int}
\let\!lim\lim \def\lim{\displaystyle\!lim}
\let\!sum\sum \def\sum{\displaystyle\!sum}
\let\!sup\sup \def\sup{\displaystyle\!sup}
\let\!inf\inf \def\inf{\displaystyle\!inf}
\let\!cap\cap \def\cap{\displaystyle\!cap}
\let\!max\max \def\max{\displaystyle\!max}
\let\!min\min \def\min{\displaystyle\!min}
\let\!frac\frac \def\frac{\displaystyle\!frac}
\catcode`!=12

\let\oldsection\section
\renewcommand\section{\setcounter{equation}{0}\oldsection}

\allowdisplaybreaks
\def\pf{\it{Proof.}\rm\quad}

\def\N{\mathbb{N}}

\newtheorem{thm}{Theorem}[section]

\newtheorem{cor}[thm]{Corollary}

\newtheorem{exa}{Example}[section]
\begin{document}
\title {\bf Evaluations of some quadratic Euler sums}
\author{
{Xin Si\thanks{Email: xsi@xmut.edu.cn (X. Si)} \quad Ce Xu\thanks{Corresponding author. Email: 15959259051@163.com (C. Xu)}}\\[1mm]
\small $\ast$School of Applied Mathematics, Xiamen University of Technology\\
\small Xiamen
361005, P.R. China\\
\small $\dagger$ School of Mathematical Sciences, Xiamen University\\
\small Xiamen
361005, P.R. China}

\date{}
\maketitle \noindent{\bf Abstract } This paper develops an approach to the evaluation of quadratic Euler sums that involve harmonic numbers. The approach is based on simple integral computations of polylogarithms.
By using the approach, we establish some relations between quadratic Euler sums and linear sums. Furthermore, we
obtain some closed form representations of quadratic sums in terms of zeta values and linear sums. The given representations are new.
\\[2mm]
\\[2mm]

\section{Introduction}
Euler sums are real numbers, originally defined by Euler sums, that have been much studied in recent years because of their many surprising properties and the many places they appear in mathematics and mathematical physics. There are many conjectures concerning the values of Euler sums, for example see [2, 6, 9, 12, 19, 20]. The subject of this paper is Euler sums, which are the infinite sums whose general term is a product of harmonic numbers (or alternating harmonic numbers) of index $n$ and a power of $n^{-1}$ (or ${(-1)^{n-1}}{n^{-1}}$).
The $n$th generalized harmonic numbers and $n$th generalized alternating harmonic numbers are defined by
$$\zeta_n(k):=\sum\limits_{j=1}^n\frac {1}{j^k} ,\ L_{n}(k):=\sum\limits_{j=1}^n\frac{(-1)^{j-1}}{j^k},\ k,n \in \N:=\{1,2,3\ldots\},\eqno(1.1)$$
where $H_n:=\zeta_n(1)=\sum\limits_{j=1}^n \frac 1{j}$ is the natural harmonic number.
The classical linear Euler sum is defined by
$$S(p;q):=\sum\limits_{n = 1}^\infty  {\frac{1}{{{n^q}}}} \sum\limits_{k = 1}^n {\frac{1}{{{k^p}}}}, \eqno(1.2)$$
where $p,q$ are positive integers with $q \geq 2$ and the quantity $w:=p+q$ is called the weight. The earliest results on linear sums $S(p;q)$ are due to Euler who elaborated a method to reduce double sums of small weight to certain rational linear combinations of products of zeta values. In particular, he proved the simple relation in 1775 (see [2, 12])
\[S(1;k)=\sum\limits_{n = 1}^\infty  {\frac{{{H_n}}}{{{n^k}}}}  = \frac{1}{2}\left\{ {\left( {k + 2} \right)\zeta \left( {k + 1} \right) - \sum\limits_{i = 1}^{k - 2} {\zeta \left( {k - i} \right)\zeta \left( {i + 1} \right)} } \right\},\eqno(1.3)\]
and determined the explicit values of zeta values function at even integers:
\[\zeta \left( {2m} \right) = \frac{{{{\left( { - 1} \right)}^m}{B_{2m}}}}{{2\left( {2m} \right)!}}{\left( {2\pi } \right)^{2m}},\]
where $B_k\in \mathbb{Q}$ are the Bernoulli numbers defined by the generating function (see [1, 4, 5])
\[\frac{x}{{{e^x} - 1}} := \sum\limits_{k = 0}^\infty  {{B_k}\frac{{{x^k}}}{{k!}}} .\]
It is easy to verify that $B_0=1,\ B_1=-\frac 1{2},\ B_2=\frac 1{6},\ B_4=-\frac 1{30},$ and $B_{2m+1}=0$ for $m\geq 1$.
The Riemann zeta function and alternating Riemann zeta function are defined respectively by
 $$\zeta(s):=\sum\limits_{n = 1}^\infty {\frac {1}{n^{s}}},\Re(s)>1$$
and
\[\bar \zeta \left( s \right) := \sum\limits_{n = 1}^\infty  {\frac{{{{\left( { - 1} \right)}^{n - 1}}}}{{{n^s}}}},{\mathop{\Re}\nolimits} \left( s \right) \ge 1.\]
Obviously, for $\Re(s)>1$, $\bar \zeta \left( s \right) = \left( {1 - \frac{1}{{{2^{s - 1}}}}} \right)\zeta \left( s \right).$
The general multiple zeta functions is defined as
\[\zeta \left( {{s_1},{s_2}, \cdots ,{s_m}} \right) := \sum\limits_{{n_1} > {n_2} >  \cdots  > {n_m} > 0}^{} {\frac{1}{{n_1^{{s_1}}n_2^{{s_2}} \cdots n_m^{{s_m}}}}}, \]
where ${s_1} +  \cdots  + {s_m}$ is called the weight and $m$ is the multiplicity.\\
Euler conjectured that the linear sums $S(p;q)$ would be reducible to zeta values whenever $p + q$ is odd, and even gave what he hoped to obtain the general formula. In \cite{BBG1995}, D. Borwein, J.M. Borwein and R. Girgensohn proved conjecture and formula,
and in [2], D.H. Bailey, J.M. Borwein and R. Girgensohn demonstrated that it is ``very likely" that linear sums with $p + q > 7,\ p + q$ even, are not reducible.\\
Next, we introduce the generalized Euler sums. For integers $q,p_1,\ldots,p_m$ with $q\geq 2$, we define the generalized Euler sums as
\[S\left( {{p_1},{p_2}, \ldots ,{p_m};q} \right): = \sum\limits_{n = 1}^\infty  {\frac{{{X_n}\left( {{p_1}} \right){X_n}\left( {{p_2}} \right) \cdots {X_n}\left( {{p_m}} \right)}}{{{n^q}}}} ,\tag{1.4}\]
where ${X_n}\left( {{p_i}} \right) = {\zeta _n}\left( {{p_i}} \right)$ if $p_i>0$, and ${X_n}\left( {{p_i}} \right) = {L_n}\left( {{-p_i}} \right)$ otherwise. In below, if $p<0$, we will denote it by $\overline {-p}$. For example,
\[S\left( {2,\bar 3,5,\bar 7;q} \right) = \sum\limits_{n = 1}^\infty  {\frac{{{\zeta _n}\left( 2 \right){L_n}\left( 3 \right){\zeta _n}\left( 5 \right){L_n}\left( 7 \right)}}{{{n^q}}}} .\]
Similarly, for $q\in \N$, we define
\[S\left( {{p_1},{p_2}, \ldots ,{p_m};\bar q} \right): = \sum\limits_{n = 1}^\infty  {\frac{{{X_n}\left( {{p_1}} \right){X_n}\left( {{p_2}} \right) \cdots {X_n}\left( {{p_m}} \right)}}{{{n^q}}}{{\left( { - 1} \right)}^{n - 1}}}.\tag{1.5} \]
We call $w: = \left| {{p_1}} \right| + \left| {{p_2}} \right| +  \cdots  + \left| {{p_m}} \right| + \left| q \right|$ the weight of the Euler sums $S(p_1,p_2,\ldots,p_m;q)$. Hence, by the definition of $S(p_1,p_2,\ldots,p_m;q)$, we know that the linear sums are altogether four types:
\begin{align*}
&S\left( {p;q} \right) = \sum\limits_{n = 1}^\infty  {\frac{{{\zeta _n}\left( p \right)}}{{{n^q}}}} ,\:S\left( {\bar p;q} \right) = \sum\limits_{n = 1}^\infty  {\frac{{{L_n}\left( p \right)}}{{{n^q}}}} ,\:\\&S\left( {p;\bar q} \right) = \sum\limits_{n = 1}^\infty  {\frac{{{\zeta _n}\left( p \right)}}{{{n^q}}}} {\left( { - 1} \right)^{n - 1}},\:S\left( {\bar p;\bar q} \right) = \sum\limits_{n = 1}^\infty  {\frac{{{L_n}\left( p \right)}}{{{n^q}}}{{\left( { - 1} \right)}^{n - 1}}}
\end{align*}
The study of these Euler sums was started by Euler. After that many different methods, including partial fraction expansions, Eulerian Beta integrals, summation formulas for generalized hypergeometric functions and contour integrals, have been used to evaluate these sums.
The relationship between the
values of the Riemann zeta function and Euler sums has been studied by many authors, for details and historical introductions, please see [2, 3, 6-20]. For example, in [12], Philippe Flajolet and Bruno Salvy gave explicit reductions to zeta values for all linear sums
$S\left( {p;q} \right),S\left( {\bar p;q} \right),S\left( {p;\bar q} \right),S\left( {\bar p;\bar q} \right)$
with $w:=p+q$ odd. Moreover, they proved the following conclusion:
If $p_1+p_2+q$ is even, and $p_1>1, p_2>1, q>1$, the quadratic sums
\[{S({{p_1},{p_2};q})} = \sum\limits_{n = 1}^\infty  {\frac{{{\zeta _n}\left( {{p_1}} \right){\zeta _n}\left( {{p_2}} \right)}}{{{n^q}}}} \]
 are reducible to linear sums (see Theorem 4.2 in the reference [12]). In [19] and [21], we showed that all quadratic sums $S(p_1,p_2;q)$ with $\left| {{p_1}} \right| + \left| {{p_2}} \right| + \left| q \right| \le 4$ were reducible to zeta values and polylogarithms,
and in [25], we proved that all Euler sums of the form $S(1,p;q)$ for weights $p+q+1\in\{4,5,6,7,9\}$ with $p\geq1$ and $ q\geq2$ are expressible polynomially in terms of zeta values. For weight 8, all such sums are the sum of a polynomial in zeta values and a rational multiple of $S(2;6)$. But so far, no one has proved whether or not the
Euler sums
\[{S({1,m;p})} = \sum\limits_{n = 1}^\infty  {\frac{{{H_n}{\zeta _n}\left( m \right)}}{{{n^p}}}} ,\;p,m \in \mathbb{N} \setminus \{1\}:=\{2,3,4,\cdots\}\ \] can be expressed in terms of zeta values and linear sums. \\
The main purpose of this paper is to evaluate some quadratic Euler sums which involve harmonic numbers and alternating harmonic numbers. In this paper, we will prove that all quadratic sums
\[S(1,p+1;p+2m)=\sum\limits_{n = 1}^\infty  {\frac{{{H_n}{\zeta _n}\left( {p + 1} \right)}}{{{n^{p + 2m}}}}},\ S(1,p+2m;p+1)= \sum\limits_{n = 1}^\infty  {\frac{{{H_n}{\zeta _n}\left( {p + 2m} \right)}}{{{n^{p + 1}}}}} ,\;p,m \in \N.\]
are reducible to polynomials in zeta values and to linear sums. Moreover, we also prove that, for $p \in \mathbb{N} \setminus \{1\},\; m \in \N\cup \{0\}$, the quadratic combinations
\begin{align*}
&S(\bar{1},p+2m+1;p)+S(\bar{1},p;p+2m+1)=\sum\limits_{n = 1}^\infty  {\left\{ {\frac{{{L_n}\left( 1 \right){\zeta _n}\left( {p + 2m + 1} \right)}}{{{n^p}}} + \frac{{{L_n}\left( 1 \right){\zeta _n}\left( p \right)}}{{{n^{p + 2m + 1}}}}} \right\}} ,\\
&S(\bar{1},p+2m;p)-S(\bar{1},p;p+2m)=\sum\limits_{n = 1}^\infty  {\left\{ {\frac{{{L_n}\left( 1 \right){\zeta _n}\left( {p + 2m} \right)}}{{{n^p}}} - \frac{{{L_n}\left( 1 \right){\zeta _n}\left( p \right)}}{{{n^{p + 2m}}}}} \right\}} \end{align*}
and
$$S(1,p+2m+2;p)-S(1,p;p+2m+2)=\sum\limits_{n = 1}^\infty  {\left\{ {\frac{{{H_n}{\zeta _n}\left( {p + 2m + 2} \right)}}{{{n^p}}} - \frac{{{H_n}{\zeta _n}\left( p \right)}}{{{n^{p + 2m + 2}}}}} \right\}} $$
reduce to linear sums and polynomials in zeta values.
\section{Main Theorems and Proofs}
In this section, by calculating the integrals of polylogarithm functions, we will establish some explicit relationships which involve quadratic sums and linear sums. The polylogarithm function is defined as follows
\[{\rm Li}{_p}\left( x \right) = \sum\limits_{n = 1}^\infty  {\frac{{{x^n}}}{{{n^p}}}}, \Re(p)>1,\ \left| x \right| \le 1 ,\]
with ${\rm Li_1}(x)=-\log(1-x),\ x\in [-1,1)$. First, we give the following Theorem, which will be useful in the development of the main results.
\begin{thm}\label{lem 2.1}
Let $m,p\ge 2$ be positive integers, the following identity holds:
\[\sum\limits_{n = 1}^\infty  {\frac{{\zeta \left( m \right){\zeta _n}\left( p \right) - \zeta \left( p \right){\zeta _n}\left( m \right)}}{n}}  = \zeta \left( p \right)\sum\limits_{n = 1}^\infty  {\frac{{{H_n}}}{{{n^m}}}}  - \zeta \left( m \right)\sum\limits_{n = 1}^\infty  {\frac{{{H_n}}}{{{n^p}}}}  + \zeta \left( m \right)\zeta \left( {p + 1} \right) - \zeta \left( p \right)\zeta \left( {m + 1} \right).\tag{2.1}\]
\end{thm}
\pf We construct the generating function
 \[y = \sum\limits_{n = 1}^\infty  {\left\{ {{H_n}{\zeta _n}\left( m \right) - {\zeta _n}\left( {m + 1} \right)} \right\}{x^{n - 1}}},\ x\in (-1,1). \tag{2.2}\]
By definition, the harmonic numbers satisfy the recurrence relation
\[{\zeta _{n + 1}}\left( m \right) = {\zeta _n}\left( m \right) + \frac{1}{{{{\left( {n + 1} \right)}^m}}}.\]
Then the sum on the right hand side of (2.2) is equal to
\begin{align*}
&\sum\limits_{n = 1}^\infty  {\left\{ {{H_n}{\zeta _n}\left( m \right) - {\zeta _n}\left( {m + 1} \right)} \right\}{x^{n - 1}}} \\&= \sum\limits_{n = 1}^\infty  {\left\{ {\left( {{H_n} + \frac{1}{{n + 1}}} \right)\left( {{\zeta _n}\left( m \right) + \frac{1}{{{{\left( {n + 1} \right)}^m}}}} \right) - \left( {{\zeta _n}\left( {m + 1} \right) + \frac{1}{{{{\left( {n + 1} \right)}^{m + 1}}}}} \right)} \right\}{x^n}}
\nonumber \\ &=\sum\limits_{n = 1}^\infty  {\left\{ {{H_n}{\zeta _n}\left( m \right) - {\zeta _n}\left( {m + 1} \right) + \frac{{{H_n}}}{{{{\left( {n + 1} \right)}^m}}} + \frac{{{\zeta _n}\left( m \right)}}{{n + 1}}} \right\}{x^n}} .
\end{align*}
By simple calculation, we get
\[\sum\limits_{n = 1}^\infty  {\left\{ {{H_n}{\zeta _n}\left( m \right) - {\zeta _n}\left( {m + 1} \right)} \right\}{x^{n - 1}}}  = \sum\limits_{n = 1}^\infty  {\left\{ {\frac{{{H_n}}}{{{{\left( {n + 1} \right)}^m}}} + \frac{{{\zeta _n}\left( m \right)}}{{n + 1}}} \right\}\frac{{{x^n}}}{{1 - x}}}.\tag{2.3} \]
Multiplying (2.3) by ${\ln ^{p - 1}}x$ and integrating over (0,1), we obtain the formula
\[\sum\limits_{n = 1}^\infty  {\frac{{{H_n}{\zeta _n}\left( m \right) - {\zeta _n}\left( {m + 1} \right)}}{{{n^p}}}}  = \sum\limits_{n = 1}^\infty  {\left\{ {\frac{{{H_n}}}{{{{\left( {n + 1} \right)}^m}}} + \frac{{{\zeta _n}\left( m \right)}}{{n + 1}}} \right\}\left\{ {\zeta \left( p \right) - {\zeta _n}\left( p \right)} \right\}}.\tag{2.4}\]
After some straightforward manipulations, formula (2.4) can be written as
\begin{align*}
\sum\limits_{n = 1}^\infty  {\left\{ {\frac{{{H_n}{\zeta _n}\left( m \right)}}{{{n^p}}} + \frac{{{H_n}{\zeta _n}\left( p \right)}}{{{n^m}}}} \right\}} &=\zeta \left( p \right)\sum\limits_{n = 1}^\infty  {\frac{{{H_n}}}{{{n^m}}}}  + \sum\limits_{n = 1}^\infty  {\frac{{{H_n}}}{{{n^{p + m}}}}}  + \sum\limits_{n = 1}^\infty  {\frac{{{\zeta _n}\left( m \right)}}{{{n^{p + 1}}}}}
\nonumber \\ &\quad- \sum\limits_{n = 1}^\infty  {\frac{{{\zeta _n}\left( {m + 1} \right)}}{{{n^p}}}}  + \sum\limits_{n = 1}^\infty  {\frac{{{\zeta _n}\left( m \right)}}{n}} \left\{ {\zeta \left( p \right) - {\zeta _n}\left( p \right)} \right\}. \tag{2.5}
\end{align*}
Change $(m,p)$ to $(p,m)$, the result is
\begin{align*}
\sum\limits_{n = 1}^\infty  {\left\{ {\frac{{{H_n}{\zeta _n}\left( p \right)}}{{{n^m}}} + \frac{{{H_n}{\zeta _n}\left( m \right)}}{{{n^p}}}} \right\}} &=\zeta \left( m \right)\sum\limits_{n = 1}^\infty  {\frac{{{H_n}}}{{{n^p}}}}  + \sum\limits_{n = 1}^\infty  {\frac{{{H_n}}}{{{n^{m + p}}}}}  + \sum\limits_{n = 1}^\infty  {\frac{{{\zeta _n}\left( p \right)}}{{{n^{m + 1}}}}}
\nonumber \\ &\quad - \sum\limits_{n = 1}^\infty  {\frac{{{\zeta _n}\left( {p + 1} \right)}}{{{n^m}}}}  + \sum\limits_{n = 1}^\infty  {\frac{{{\zeta _n}\left( p \right)}}{n}} \left\{ {\zeta \left( m \right) - {\zeta _n}\left( m \right)} \right\}. \tag{2.6}
\end{align*}
Therefore, combining (2.5) and (2.6), we obtain the desired result. The proof of Theorem 2.1 is completed.  \hfill$\square$\\
Proceeding in a similar fashion to evaluation of the Theorem 2.1, we consider the following function
\[y = \sum\limits_{n = 1}^\infty  {\left\{ {{\zeta _n}\left( {1,a+1} \right){\zeta _n}\left( {p,a+1} \right) - {\zeta _n}\left( {p + 1,a+1} \right)} \right\}{x^{n + a - 1}}} ,\;x \in \left( { - 1,1} \right)\;\]
where the partial sums ${{\zeta _n}\left( {p,a+1} \right)}$ for $p\geq1$ of Hurwitz zeta function is defined as
\[{\zeta _n}\left( {p,a+1} \right) := \sum\limits_{k = 1}^n {\frac{1}{{{{\left( {k + a} \right)}^p}}}} ,\;a \notin \N^-:= \left\{ {-1,-2,-3 \cdots } \right\}.\]
The Hurwitz zeta function is defined by
\[\zeta \left( {p,a + 1} \right) := \sum\limits_{n = 1}^\infty  {\frac{1}{{{{\left( {n + a} \right)}^p}}}} ,\Re\left( p \right) > 1,a \notin  \N^-.\]
By a similar argument as in the proof of Theorem 2.1, we deduce the more general identity
\begin{align*}
&\sum\limits_{n = 1}^\infty  {\frac{{\zeta \left( {m,a + 1} \right){\zeta _n}\left( {p,a+1} \right) - \zeta \left( {p,a + 1} \right){\zeta _n}\left( {m,a+1} \right)}}{{n + a}}} \\
& = \zeta \left( {p,a + 1} \right)\sum\limits_{n = 1}^\infty  {\frac{{{\zeta _n}\left( {1,a+1} \right)}}{{{{\left( {n + a} \right)}^m}}}}  - \zeta \left( {m,a + 1} \right)\sum\limits_{n = 1}^\infty  {\frac{{{\zeta _n}\left( {1,a+1} \right)}}{{{{\left( {n + a} \right)}^p}}}} \\
&\quad + \zeta \left( {m,a + 1} \right)\zeta \left( {p + 1,a + 1} \right) - \zeta \left( {m + 1,a + 1} \right)\zeta \left( {p,a + 1} \right).
\end{align*}
When $a=0$, the result is formula (2.1).
\begin{thm}\label{lem 2.2}\ \  Let $p \geq 2,m \ge 0$ be integers and $x \in \left[ { - 1,1} \right)$. Then the following identity holds:
\begin{align*}
&{\left( { - 1} \right)^{p - 1}}\sum\limits_{n = 1}^\infty  {\left\{ {\frac{{{\zeta _n}\left( {p + 2m + 1} \right)}}{{{n^p}}} + \frac{{{\zeta _n}\left( p \right)}}{{{n^{p + 2m + 1}}}}} \right\}\left( {\sum\limits_{k = 1}^n {\frac{{{x^k}}}{k}} } \right)}\\  &= \sum\limits_{i = 1}^{p + 2m} {{{\left( { - 1} \right)}^{i - 1}}{\rm{L}}{{\rm{i}}_{p + 2m + 2 - i}}\left( x \right)} \sum\limits_{n = 1}^\infty  {\frac{{{\zeta _n}\left( p \right)}}{{{n^i}}}} {x^n} - \sum\limits_{i = 1}^{p - 1} {{{\left( { - 1} \right)}^{i - 1}}{\rm{L}}{{\rm{i}}_{p + 1 - i}}\left( x \right)} \sum\limits_{n = 1}^\infty  {\frac{{{\zeta _n}\left( {p + 2m + 1} \right)}}{{{n^i}}}} {x^n}
\nonumber \\ &\quad  + {\left( { - 1} \right)^p}\ln \left( {1 - x} \right)\sum\limits_{n = 1}^\infty  {\left\{ {\frac{{{\zeta _n}\left( {p + 2m + 1} \right)}}{{{n^p}}} + \frac{{{\zeta _n}\left( p \right)}}{{{n^{p + 2m + 1}}}}} \right\}\left( {1 - {x^n}} \right)} . \tag{2.7}
\end{align*}
\end{thm}
\pf By the definition of polylogarithm function and Cauchy product formula, we can verify that
\[\frac{{{\rm Li}{_m}\left( x \right)}}{{1 - x}} = \sum\limits_{n = 1}^\infty  {{\zeta _n}\left( m \right){x^n}}, \ x\in (-1,1).\tag{2.8}\]
Now, we consider the integral \[\int\limits_0^x {\frac{{{\rm{L}}{{\rm{i}}_p}\left( t \right){\rm{L}}{{\rm{i}}_{p + 2m + 1}}\left( t \right)}}{{t\left( {1 - t} \right)}}} dt,\;x \in \left( { - 1,1} \right).\]
First, by virtue of (2.8), we obtain
\[\int\limits_0^x {\frac{{{\rm{L}}{{\rm{i}}_p}\left( t \right){\rm{L}}{{\rm{i}}_{p + 2m + 1}}\left( t \right)}}{{t\left( {1 - t} \right)}}} dt = \sum\limits_{n = 1}^\infty  {{\zeta _n}\left( p \right)} \int\limits_0^x {{t^{n - 1}}{\rm{L}}{{\rm{i}}_{p + 2m + 1}}\left( t \right)} dt = \sum\limits_{n = 1}^\infty  {{\zeta _n}\left( {p + 2m + 1} \right)} \int\limits_0^x {{t^{n - 1}}{\rm{L}}{{\rm{i}}_p}\left( t \right)} dt.\tag{2.9}\]
On the other hand, using integration by parts we deduce that
\[\int\limits_0^x {{t^{n - 1}}{\rm Li}{_p}\left( t \right)dt}  = \sum\limits_{i = 1}^{p - 1} {{{\left( { - 1} \right)}^{i - 1}}\frac{{{x^n}}}{{{n^i}}}{\rm Li}{_{p + 1 - i}}\left( x \right)}  + \frac{{{{\left( { - 1} \right)}^p}}}{{{n^p}}}\ln \left( {1 - x} \right)\left( {{x^n} - 1} \right) - \frac{{{{\left( { - 1} \right)}^p}}}{{{n^p}}}\left( {\sum\limits_{k = 1}^n {\frac{{{x^k}}}{k}} } \right).\tag{2.10}\]
In fact, by using the elementary integral identity
\[\int\limits_0^1 {{x^{n - 1}}{{\ln }^m}x\ln \left( {1 - x} \right)} dx = {\left( { - 1} \right)^{m + 1}}m!\left\{ {\frac{{{H_n}}}{{{n^{m + 1}}}} - \sum\limits_{j = 1}^m {\frac{{\zeta \left( {j + 1} \right) - {\zeta _n}\left( {j + 1} \right)}}{{{n^{m + 1 - j}}}}} } \right\},\]
then multiplying (2.10) by $\frac{{{{\ln }^{m-1}}x}}{{ x}}$ and integrating over the interval $(0,1)$, we have the following recurrence relation
\begin{align*}
\int\limits_0^1 {{x^{n - 1}}{{\ln }^m}x{\rm{L}}{{\rm{i}}_p}\left( x \right)} dx =& m\sum\limits_{i = 1}^{p - 1} {\frac{{{{\left( { - 1} \right)}^i}}}{{{n^i}}}} \int\limits_0^1 {{x^{n - 1}}{{\ln }^{m - 1}}x{\rm{L}}{{\rm{i}}_{p + 1 - i}}\left( x \right)} dx + m!{\left( { - 1} \right)^{m + p - 1}}\frac{{{\zeta _n}\left( {m + 1} \right)}}{{{n^p}}}\\
&\quad + m!\frac{{{{\left( { - 1} \right)}^{m + p - 1}}}}{{{n^p}}}\left\{ {\frac{{{H_n}}}{{{n^m}}} - \sum\limits_{j = 1}^{m - 1} {\frac{{\zeta \left( {j + 1} \right) - {\zeta _n}\left( {j + 1} \right)}}{{{n^{m - j}}}}}  - \zeta \left( {m + 1} \right)} \right\}.\tag{2.11}
\end{align*}
Substituting (2.10) into (2.9), we get
\begin{align*}
&\int\limits_0^x {\frac{{{\rm{L}}{{\rm{i}}_p}\left( t \right){\rm{L}}{{\rm{i}}_{p + 2m + 1}}\left( t \right)}}{{t\left( {1 - t} \right)}}} dt\\
& =\sum\limits_{i = 1}^{p + 2m} {{{\left( { - 1} \right)}^{i - 1}}{\rm{L}}{{\rm{i}}_{p + 2m + 2 - i}}\left( x \right)} \sum\limits_{n = 1}^\infty  {\frac{{{\zeta _n}\left( p \right)}}{{{n^i}}}} {x^n} + {\left( { - 1} \right)^p}\ln \left( {1 - x} \right)\sum\limits_{n = 1}^\infty  {\frac{{{\zeta _n}\left( p \right)}}{{{n^{p + 2m + 1}}}}\left( {1 - {x^n}} \right)} \\
&\quad + {\left( { - 1} \right)^p}\sum\limits_{n = 1}^\infty  {\frac{{{\zeta _n}\left( p \right)}}{{{n^{p + 2m + 1}}}}\left( {\sum\limits_{k = 1}^n {\frac{{{x^k}}}{k}} } \right)} \\
&=\sum\limits_{i = 1}^{p - 1} {{{\left( { - 1} \right)}^{i - 1}}{\rm{L}}{{\rm{i}}_{p + 1 - i}}\left( x \right)} \sum\limits_{n = 1}^\infty  {\frac{{{\zeta _n}\left( {p + 2m + 1} \right)}}{{{n^i}}}} {x^n} + {\left( { - 1} \right)^{p - 1}}\ln \left( {1 - x} \right)\sum\limits_{n = 1}^\infty  {\frac{{{\zeta _n}\left( {p + 2m + 1} \right)}}{{{n^p}}}\left( {1 - {x^n}} \right)}   \\
&\quad  + {\left( { - 1} \right)^{p - 1}}\sum\limits_{n = 1}^\infty  {\frac{{{\zeta _n}\left( {p + 2m + 1} \right)}}{{{n^p}}}\left( {\sum\limits_{k = 1}^n {\frac{{{x^k}}}{k}} } \right)} .\tag{2.12}
\end{align*}
By a direct calculation, we deduce the result.\hfill$\square$\\
Noting that when $x$ approach 1, by using (2.1), we obtain the result
\begin{align*}
&\mathop {\lim }\limits_{x \to 1} \left\{ {{\rm{L}}{{\rm{i}}_m}\left( x \right)\sum\limits_{n = 1}^\infty  {\frac{{{\zeta _n}\left( p \right)}}{n}{x^n}}  - {\rm{L}}{{\rm{i}}_p}\left( x \right)\sum\limits_{n = 1}^\infty  {\frac{{{\zeta _n}\left( m \right)}}{n}{x^n}} } \right\} \\
&= \sum\limits_{n = 1}^\infty  {\frac{{\zeta \left( m \right){\zeta _n}\left( p \right) - \zeta \left( p \right){\zeta _n}\left( m \right)}}{n}}\\
&=\zeta \left( p \right)S\left( {1;m} \right) - \zeta \left( m \right)S\left( {1;p} \right) + \zeta \left( m \right)\zeta \left( {p + 1} \right) - \zeta \left( p \right)\zeta \left( {m + 1} \right),\tag{2.13}
\end{align*}
where $p,m\in \mathbb{N} \setminus \{1\}.$
Hence, letting $x\rightarrow1$ and $x\rightarrow-1$ in Theorem 2.2 and combining (2.13), we get the following results.
\begin{cor}
Let $p \geq 2,m \ge 0$  be integers, we have
\begin{align*}
&{(-1)^{p-1}}\{S\left( {1,p + 2m + 1;p} \right) + S\left( {1,p;p + 2m + 1} \right)\} \\
&=\zeta \left( p \right)S(1;p+2m+1) - \zeta \left( {p + 2m + 1} \right)S(1;p) \\
&\quad+ \zeta \left( {p + 1} \right)\zeta \left( {p + 2m + 1} \right) - \zeta \left( p \right)\zeta \left( {p + 2m + 2} \right)
 \\ &\quad + \sum\limits_{i = 2}^{p + 2m} {{{\left( { - 1} \right)}^{i - 1}}\zeta \left( {p + 2m + 2 - i} \right)}S(p;i)\\
 &\quad- \sum\limits_{i = 2}^{p - 1} {{{\left( { - 1} \right)}^{i - 1}}\zeta \left( {p + 1 - i} \right)} S(p+2m+1;i). \tag{2.14}
\end{align*}
\end{cor}
\begin{cor} Let $p \geq 2,m \ge 0$  be integers, we have
\begin{align*}
 &{\left( { - 1} \right)^p}\left\{ {S\left( {\bar 1,p + 2m+1;p} \right)+ S\left( {\bar 1,p;p + 2m+1} \right)} \right\}\\
 &=\sum\limits_{i = 1}^{p + 2m} {{{\left( { - 1} \right)}^{i - 1}}\bar \zeta \left( {p + 2m + 2 - i} \right)}S(p;{\bar i}) \\
           &\quad  - \sum\limits_{i = 1}^{p - 1} {{{\left( { - 1} \right)}^{i - 1}}\bar \zeta \left( {p + 1 - i} \right)}S(p+2m+1;{\bar i}) \\
  &\quad+ {\left( { - 1} \right)^p}\ln 2 \{S(p+2m+1;p)+S(p;p+2m+1)\}\\
  &\quad+ {\left( { - 1} \right)^p}\ln 2 \{S(p+2m+1;{\bar p})+S(p;\overline {p+2m+1})\}.\tag{2.15}
\end{align*}
\end{cor}
In the same way as in the proof of (2.7), we obtain the following Theorem.
\begin{thm}
For $p\in \mathbb{N} \setminus \{1\},\ m\in \N\cup \{0\}$ and $x \in \left[ { - 1,1} \right)$
. Then the following identity holds:
\begin{align*}
&{\left( { - 1} \right)^{p - 1}}\sum\limits_{n = 1}^\infty  {\left\{ {\frac{{{\zeta _n}\left( {p + 2m} \right)}}{{{n^p}}} - \frac{{{\zeta _n}\left( p \right)}}{{{n^{p + 2m}}}}} \right\}\left( {\sum\limits_{k = 1}^n {\frac{{{x^k}}}{k}} } \right)}  \\
& =\sum\limits_{i = 1}^{p + 2m - 1} {{{\left( { - 1} \right)}^{i - 1}}{\rm{L}}{{\rm{i}}_{p + 2m + 1 - i}}\left( x \right)} \sum\limits_{n = 1}^\infty  {\frac{{{\zeta _n}\left( p \right)}}{{{n^i}}}} {x^n} - \sum\limits_{i = 1}^{p - 1} {{{\left( { - 1} \right)}^{i - 1}}{\rm{L}}{{\rm{i}}_{p + 1 - i}}\left( x \right)} \sum\limits_{n = 1}^\infty  {\frac{{{\zeta _n}\left( {p + 2m} \right)}}{{{n^i}}}} {x^n}\\
&\quad + {\left( { - 1} \right)^p}\ln \left( {1 - x} \right)\sum\limits_{n = 1}^\infty  {\left\{ {\frac{{{\zeta _n}\left( {p + 2m} \right)}}{{{n^p}}} - \frac{{{\zeta _n}\left( p \right)}}{{{n^{p + 2m}}}}} \right\}\left( {1 - {x^n}} \right)}  .\tag{2.16}
\end{align*}
\end{thm}
\pf Similarly to the proof of Theorem 2.2, considering integral
\[\int\limits_0^x {\frac{{{\rm{L}}{{\rm{i}}_p}\left( t \right){\rm{L}}{{\rm{i}}_{p + 2m}}\left( t \right)}}{{t\left( {1 - t} \right)}}} dt,\:x \in \left( { - 1,1} \right).\]
Then with the help of formula (2.10) we may easily deduce the result.\hfill$\square$\\
Similarly, in (2.16), taking $x\rightarrow 1$ and $x \rightarrow -1$, by using (2.13), we can give the following Corollaries.
\begin{cor}
For integers $p\in \mathbb{N} \setminus \{1\}$ and $m\in \N\cup\{0\}$, we have
\begin{align*}
&{\left( { - 1} \right)^{p - 1}}\left\{ {S\left( {1,p + 2m;p} \right) - S\left( {1,p;p + 2m} \right)} \right\} \\
& =\zeta \left( p \right)S(1;p+2m) - \zeta \left( {p + 2m} \right)S(1;p) \\
&\quad  + \zeta \left( {p + 1} \right)\zeta \left( {p + 2m} \right) - \zeta \left( p \right)\zeta \left( {p + 2m + 1} \right)\\
&\quad+\sum\limits_{i = 2}^{p + 2m - 1} {{{\left( { - 1} \right)}^{i - 1}}\zeta \left( {p + 2m + 1 - i} \right)} S(p;i)  \\&\quad- \sum\limits_{i = 2}^{p - 1} {{{\left( { - 1} \right)}^{i - 1}}\zeta \left( {p + 1 - i} \right)} S(p+2m;i)
.\tag{2.17}
\end{align*}
\end{cor}
\begin{cor}
For integers $p\in \mathbb{N} \setminus \{1\}$ and $m\in \N\cup\{0\}$, we have
\begin{align*}
&{\left( { - 1} \right)^p}\left\{ {S\left( {\bar 1,p + 2m;p} \right) - S\left( {\bar 1,p;p + 2m} \right)} \right\} \\
& =\sum\limits_{i = 1}^{p + 2m - 1} {{{\left( { - 1} \right)}^{i - 1}}\bar \zeta \left( {p + 2m + 1 - i} \right)}S(p;\bar i)  \\
&\quad - \sum\limits_{i = 1}^{p - 1} {{{\left( { - 1} \right)}^{i - 1}}\bar \zeta \left( {p + 1 - i} \right)}S(p+2m;\bar i) \\
&\quad + {\left( { - 1} \right)^p}\ln 2\left\{ {S\left( {p + 2m;p} \right) - S\left( {p;p + 2m} \right) } \right\}\\
&\quad + {\left( { - 1} \right)^p}\ln 2\left\{ { S\left( {p + 2m;\bar p} \right) - S\left( {p;\overline {p + 2m}} \right)} \right\}
 .  \tag{2.18}
\end{align*}
\end{cor}
\begin{thm}
For $ {l_1},{l_2},m \in \N$ and $x,y,z \in \left[ { - 1,1} \right)$, we have the following relation
\begin{align*}
&\sum\limits_{n = 1}^\infty  {\frac{{{\zeta _n}\left( {{l_1};x} \right){\zeta _n}\left( {{l_2};y} \right)}}{{{n^m}}}{z^n}}  + \sum\limits_{n = 1}^\infty  {\frac{{{\zeta _n}\left( {{l_1};x} \right){\zeta _n}\left( {m;z} \right)}}{{{n^{{l_2}}}}}{y^n}}  + \sum\limits_{n = 1}^\infty  {\frac{{{\zeta _n}\left( {{l_2};y} \right){\zeta _n}\left( {m;z} \right)}}{{{n^{{l_1}}}}}{x^n}} \\
& =\sum\limits_{n = 1}^\infty  {\frac{{{\zeta _n}\left( {m;z} \right)}}{{{n^{{l_1} + {l_2}}}}}{{\left( {xy} \right)}^n}}  + \sum\limits_{n = 1}^\infty  {\frac{{{\zeta _n}\left( {{l_1};x} \right)}}{{{n^{m + {l_2}}}}}{{\left( {yz} \right)}^n}}  + \sum\limits_{n = 1}^\infty  {\frac{{{\zeta _n}\left( {{l_2};y} \right)}}{{{n^{{l_1} + m}}}}{{\left( {xz} \right)}^n}}  \\
&\quad  +{\rm Li}{_m}\left( z \right){\rm Li}{_{{l_1}}}\left( x \right){\rm Li}{_{{l_2}}}\left( y \right) - {\rm Li}{_{{l_1} + {l_2} + m}}\left( {xyz} \right) \tag{2.19}
\end{align*}
where the partial sum ${\zeta _n}\left( {l;x} \right)$ is defined by ${\zeta _n}\left( {l;x} \right) := \sum\limits_{k = 1}^n {\frac{{{x^k}}}{k^l}}$.
\end{thm}
\pf  We construct the function $F\left( {x,y,z} \right) = \sum\limits_{n = 1}^\infty  {\left\{ {{\zeta _n}\left( {{l_1};x} \right){\zeta _n}\left( {{l_2};y} \right) - {\zeta _n}\left( {{l_1} + {l_2};xy} \right)} \right\}{z^{n - 1}}} $. By the definition of ${\zeta _n}\left( {l;x} \right)$, we have
\[F\left( {x,y,z} \right) = zF\left( {x,y,z} \right) + \sum\limits_{n = 1}^\infty  {\left\{ {\frac{{{\zeta _n}\left( {{l_1};x} \right)}}{{{{\left( {n + 1} \right)}^{{l_2}}}}}{y^{n+1}} + \frac{{{\zeta _n}\left( {{l_2};y} \right)}}{{{{\left( {n + 1} \right)}^{{l_1}}}}}{x^{n+1}}} \right\}{z^n}} .\tag{2.20}\]
Moving $zF(x,y,z)$ from right to left and then multiplying  $(1-z)^{-1}$  to the equation (2.20) and integrating over the interval $(0,z)$, we obtain
\begin{align*}
&\sum\limits_{n = 1}^\infty  {\frac{{{\zeta _n}\left( {{l_1};x} \right){\zeta _n}\left( {{l_2};y} \right) - {\zeta _n}\left( {{l_1} + {l_2};xy} \right)}}{n}{z^{n}}} \\ &= \sum\limits_{n = 1}^\infty  {\left\{ {\frac{{{\zeta _n}\left( {{l_1};x} \right)}}{{{{\left( {n + 1} \right)}^{{l_2}}}}}{y^{n+1}} + \frac{{{\zeta _n}\left( {{l_2};y} \right)}}{{{{\left( {n + 1} \right)}^{{l_1}}}}}{x^{n+1}}} \right\}\left\{ {{\rm Li}{_1}\left( z \right) - {\zeta _n}\left( {1;z} \right)} \right\}}. \tag{2.21}
\end{align*}
Furthermore, using integration and the following formula
\[\sum\limits_{n = 1}^\infty  {\left\{ {\frac{{{\zeta _n}\left( {{l_1};x} \right)}}{{{{\left( {n + 1} \right)}^{{l_2}}}}}{y^{n + 1}} + \frac{{{\zeta _n}\left( {{l_2};y} \right)}}{{{{\left( {n + 1} \right)}^{{l_1}}}}}{x^{n + 1}}} \right\}}  = {\rm Li}{_{{l_1}}}\left( x \right){\rm Li}{_{{l_2}}}\left( y \right) - {\rm Li}{_{{l_1} + {l_2}}}\left( {xy} \right),\]
we can obtain (2.19). \hfill$\square$
\\Putting $(x,y,z)=(-1,1,1),\ (l_1,l_2,m)=(1,p+2m+1,p)$ and $(x,y,z)=(-1,-1,-1),\ (l_1,l_2,m)=(1,p+2m+1,p)$ in (2.19), we can give the following Corollaries.
\begin{cor}
For integers $p\in \mathbb{N} \setminus \{1\}$ and $m\in \N\cup \{0\}$, then the following identity holds:
\begin{align*}
&S\left( {\bar 1,p + 2m + 1;p} \right) + S\left( {\bar 1,p;p + 2m + 1} \right) + S\left( {p,p + 2m + 1;\bar 1} \right) \\
& =S(p;\overline {p+2m+2})+S(\bar 1;2p+2m+1)+S(p+2m+1;\overline {p+1}) \\
&\quad  +\ln 2\zeta \left( {p + 2m + 1} \right)\zeta \left( p \right) - \bar \zeta \left( {2p + 2m + 2} \right). \tag{2.22}
\end{align*}
\end{cor}
\begin{cor}
For integers $ p\in \N$ and $m\in \N\cup \{0\}$, then the following identity holds:
\begin{align*}
&S(\bar 1,\overline {p+2m+1};\bar p)+S(\bar 1,\overline {p};\overline {p+2m+1})+S(\bar p,\overline {p+2m+1};\bar 1)  \\
& =S(\bar p; p+2m+2)+S(\bar 1; 2p+2m+1)+S(\overline {p+2m+1};p+1) \\
&\quad  +\ln 2\bar \zeta \left( {p + 2m + 1} \right)\bar \zeta \left( p \right) - \bar \zeta \left( {2p + 2m + 2} \right). \tag{2.23}
\end{align*}
\end{cor}
\section{Closed form of Euler sums}
In this section, we give some linear relations among quadratic Euler sums by using Theorem 2.2 and Theorem 2.5.
We now give the following theorems.
\begin{thm} For integers $p\in \mathbb{N} \setminus \{1\}$ and $ m \in \N\cup \{0\}$, we have
\begin{align*}
&{\left( { - 1} \right)^{p - 1}}\left\{ {S\left( {2,p + 2m + 1;p} \right) + S\left( {2,p;p + 2m + 1} \right)} \right\}\\
& =\sum\limits_{i = 1}^{p + 2m} {\sum\limits_{j = 1}^{p + 2m + 1 - i} {{{\left( { - 1} \right)}^{i + j}}\zeta \left( {p + 2m + 3 - i - j} \right)S(p;i+j) } } \\
&\quad - \sum\limits_{i = 1}^{p - 1} {\sum\limits_{j = 1}^{p - i} {{{\left( { - 1} \right)}^{i + j}}\zeta \left( {p + 2 - i - j} \right)S(p+2m+1;i+j) } }  \\
&\quad - {\left( { - 1} \right)^p}\zeta \left( 2 \right)\left\{ {\zeta \left( p \right)\zeta \left( {p + 2m + 1} \right) + \zeta \left( {2p + 2m + 1} \right)} \right\}  \\
&\quad + {\left( { - 1} \right)^p}\left( {p + 2m + 1} \right)S(1,p;p+2m+2) + {\left( { - 1} \right)^p}pS(1,p+2m+1;p+1) .\tag{3.1}
\end{align*}
\end{thm}
\pf Multiplying (2.7) by $\frac 1{x}$ and integrating over (0,1), and using (2.10), we deduce Theorem 3.1 holds. \hfill$\square$
\begin{thm} For integers $p\in \mathbb{N} \setminus \{1\}$ and $ m \in \N\cup \{0\}$, we have
\begin{align*}
&{\left( { - 1} \right)^{p - 1}}\left\{ {S\left( {2,p + 2m;p} \right) - S\left( {2,p;p + 2m} \right)} \right\}\\
& =\sum\limits_{i = 1}^{p + 2m - 1} {\sum\limits_{j = 1}^{p + 2m - i} {{{\left( { - 1} \right)}^{i + j}}\zeta \left( {p + 2m + 2 - i - j} \right)S(p;i+j) } }  \\
&\quad - \sum\limits_{i = 1}^{p - 1} {\sum\limits_{j = 1}^{p - i} {{{\left( { - 1} \right)}^{i + j}}\zeta \left( {p + 2 - i - j} \right)S(p+2m;i+j) } }   \\
&\quad - {\left( { - 1} \right)^{p - 1}}\zeta \left( 2 \right)\{S(p;p+2m)-S(p+2m;p)\}  \\
&\quad + {\left( { - 1} \right)^{p - 1}}\left( {p + 2m} \right)S(1,p;p+2m+1) - {\left( { - 1} \right)^{p - 1}}pS(1,p+2m;p+1).\tag{3.2}
\end{align*}
\end{thm}
\pf By a similar argument as in the proof of Theorem 3.1, multiplying (2.16) by $\frac 1{x}$ and integrating over (0,1), and combining (2.10), we deduce Theorem 3.2 holds.\hfill$\square$\\
From [12, 20], we know that for $p\in \mathbb{N} \setminus \{1\}$ and $ m \in \N\cup \{0\}$, the quadratic sums
\begin{align*}
S(1,2;2m+1)=\frac{{{H_n}{\zeta _n}\left( 2 \right)}}{{{n^{2m+1}}}},S(2,p+2m;p)=\frac{{{\zeta _n}\left( 2 \right){\zeta _n}\left( {p + 2m} \right)}}{{{n^p}}},S(2,p:p+2m)=\frac{{{\zeta _n}\left( 2 \right){\zeta _n}\left( p \right)}}{{{n^{p + 2m}}}},
\end{align*}
are reducible to linear sums. Hence, from (3.2), we have the corollary.
\begin{cor}
For integers $p\in \mathbb{N} \setminus \{1\}$ and $ m \in \N\cup \{0\}$, the quadratic combination\[(p+2m)S(1,p;p+2m+1)-pS(1,p+2m;p+1)\] are reducible to linear sums and to polynomials in zeta values.
 \end{cor}
 On the other hand, in Corollary 2.3, we prove that for integers $p\in \mathbb{N} \setminus \{1\}$ and $ m \in \N\cup \{0\}$, the quadratic combination
\[S(1,p+2m+1;p)+S(1,p;p+2m+1)\]
can be expressed as a rational linear combination of products of zeta values and linear sums. Replacing $p$ by $p+1$ and $m$ by $m-1$ in Corollary 2.3, we obtain the following corollary.
\begin{cor}
For $p,m\in \N$, the combination
\[S(1,p+2m;p+1)+S(1,p+1;p+2m)\]
is a rational linear combination of products of zeta values and linear sums.
\end{cor}
Moreover, we note that
\begin{align*}
&\left( {p + 2m} \right)S\left( {1,p;p + 2m + 1} \right) + pS\left( {1,p + 1;p + 2m} \right)\\
=&\left\{ {\left( {p + 2m} \right)S\left( {1,p;p + 2m + 1} \right) - pS\left( {1,p + 2m;p + 1} \right)} \right\}\\
& + p\left\{ {S\left( {1,p + 2m;p + 1} \right) + S\left( {1,p + 1;p + 2m} \right)} \right\}.
\end{align*}
Therefore, from Corollary 3.3 and Corollary 3.4, we know that the combination
\[\left( {p + 2m} \right)S\left( {1,p;p + 2m + 1} \right) + pS\left( {1,p + 1;p + 2m} \right)\]
are reducible to linear sums with $p\in \mathbb{N} \setminus \{1\}$ and $m\in \N$. Since the quadratic sums $S(1,2;2m+1)$ reduce to linear sums and polynomials in zeta values. So, we obtain the following description of quadratic Euler sums $S(1,p+1;p+2m)$ and $S(1,p+2m;p+1)$.
\begin{thm} For integers $p \in \N$ and $m \in \N$, the quadratic sums
\[S(1,p+1;p+2m)=\sum\limits_{n = 1}^\infty  {\frac{{{H_n}{\zeta _n}\left( {p + 1} \right)}}{{{n^{p + 2m}}}}} ,\;S(1,p+2m;p+1)=\sum\limits_{n = 1}^\infty  {\frac{{{H_n}{\zeta _n}\left( {p + 2m} \right)}}{{{n^{p + 1}}}}} \]
are reducible to linear sums.
\end{thm}
In the following examples we collect the high-order results of quadratic Euler sums. We used the
following identities which can be easily derived from (2.14) and (3.2).
\begin{exa}
Some illustrative examples follow.
\begin{align*}
&S(1,2;3)=-\frac{101}{48}\zeta\left( 6 \right)+\frac{5}{2}{\zeta ^2}\left( 3 \right),\\
&S(1,3;2)= \frac{{227}}{{48}}\zeta \left( 6 \right) - \frac{3}{2}{\zeta ^2}\left( 3 \right),\\
&S(1,2;5)=  - \frac{{343}}{{48}}\zeta \left( 8 \right) + 12\zeta \left( 3 \right)\zeta \left( 5 \right) - \frac{5}{2}\zeta \left( 2 \right){\zeta ^2}\left( 3 \right) - \frac{3}{4}S({2;6}),\\
&S(1,3;4) =  - \frac{{511}}{{144}}\zeta \left( 8 \right) + 7\zeta \left( 3 \right)\zeta \left( 5 \right) + \zeta \left( 2 \right){\zeta ^2}\left( 3 \right) - \frac{{25}}{4}S({2;6}),\\
&S(1,4;3)= \frac{{443}}{{48}}\zeta \left( 8 \right) - \frac{{21}}{2}\zeta \left( 3 \right)\zeta \left( 5 \right) - \frac{1}{2}\zeta \left( 2 \right){\zeta ^2}\left( 3 \right) + \frac{{25}}{4}S({2;6}),\\
&S(1,5;2)=  \frac{{1063}}{{144}}\zeta \left( 8 \right) - \frac{{13}}{2}\zeta \left( 3 \right)\zeta \left( 5 \right) + \zeta \left( 2 \right){\zeta ^2}\left( 3 \right) + \frac{3}{4}S({2;6}),\\
&S(1,2;7)=  - \frac{{1331}}{{80}}\zeta( {10}) + \frac{{43}}{4}{\zeta ^2}( 5) + \frac{{41}}{2}\zeta( 3)\zeta( 7)
- 7\zeta( 2)\zeta( 3)\zeta( 5)- 2{\zeta ^2}( 3)\zeta( 4)
- \frac{5}{4}{S({2;8})},\\
&S(1,3;6) = - \frac{{247}}
{{40}}\zeta( {10}) - \frac{5}
{4}{\zeta ^2}( 5) - \frac{{15}}
{2}\zeta( 3)\zeta( 7)+ 12\zeta( 2)\zeta( 3)\zeta( 5)- \frac{{21}}
{4}{S({2,8})} - \zeta( 2){S({2,6})},\\
&S(1,4;5) = \frac{{6033}}
{{160}}\zeta( {10}) - 14{\zeta ^2}( 5) - 4\zeta( 3)\zeta( 7) - 15\zeta( 2 )\zeta( 3)\zeta( 5)- \frac{1}
{2}{\zeta ^2}( 3)\zeta( 4) + \frac{{21}}
{2}{S({2;8})} + \frac{5}
{2}\zeta( 2){S({2;6})},\\
&S(1,5;4) =  - \frac{{6569}}
{{240}}\zeta( {10}) + 16{\zeta ^2}( 5) + 10\zeta( 3)\zeta( 7) + 4\zeta( 2)\zeta( 3)\zeta( 5)+ {\zeta ^2}( 3)\zeta( 4) - \frac{{21}}
{2}{S({2;8})},\\
&S(1,6;3)  = \frac{{1043}}
{{160}}\zeta( {10}) - \frac{{17}}
{4}{\zeta ^2}( 5) - \frac{{15}}
{2}\zeta( 3)\zeta( 7) + 4\zeta( 2)\zeta( 3)\zeta( 5)- \frac{1}
{2}{\zeta ^2}( 3)\zeta( 4) - \frac{5}
{2}\zeta( 2){S({2;6})} + \frac{{21}}
{4}{S({2;8})},\\
&S(1,7;2) = \frac{{242}}
{{15}}\zeta( {10}) - \frac{{25}}
{4}{\zeta ^2}( 5) - \frac{{19}}
{2}\zeta( 3)\zeta( 7)+{\zeta ^2}( 3)\zeta( 4)+\zeta( 2){S({2;6})} + \frac{5}
{4}{S({2;8})}.
\end{align*}
\end{exa}
From (2.15) and (2.22), we have the following corollary.\\
\begin{cor} For $p\in \mathbb{N} \setminus \{1\}$ and $m \in \N\cup \{0\}$, the alternating quadratic sums
\[S(p,p+2m+1;\bar 1)=\sum\limits_{n = 1}^\infty  {\frac{{{\zeta _n}\left( p \right){\zeta _n}\left( {p + 2m + 1} \right)}}{n}} {\left( { - 1} \right)^{n - 1}}\]
are reducible to linear sums. We have
\begin{align*}
S(p,p+2m+1;\bar 1)=& S(p;\overline {p+2m+2})+S(\bar 1;2p+2m+1)+S(p+2m+1;\overline {p+1}) \\
&\quad  +\ln 2\zeta \left( {p + 2m + 1} \right)\zeta \left( p \right) - \bar \zeta \left( {2p + 2m + 2} \right)\\
&\quad+{(-1)^{p-1}}\sum\limits_{i = 1}^{p + 2m} {{{\left( { - 1} \right)}^{i - 1}}\bar \zeta \left( {p + 2m + 2 - i} \right)}S(p;{\bar i}) \\
           &\quad  -{(-1)^{p-1}} \sum\limits_{i = 1}^{p - 1} {{{\left( { - 1} \right)}^{i - 1}}\bar \zeta \left( {p + 1 - i} \right)}S(p+2m+1;{\bar i}) \\
  &\quad-\ln 2 \{S(p+2m+1;p)+S(p;p+2m+1)\}\\
  &\quad-\ln 2 \{S(p+2m+1;{\bar p})+S(p;\overline {p+2m+1})\}.\tag{3.3}
\end{align*}
\end{cor}
Letting $p=2,\ m=0$ in (2.15) and (3.3), we obtain
\begin{align*}
&S(\bar 1,3;2)+S(\bar 1,2;3)=\frac{3}{4}{\zeta ^2}\left( 3 \right) + \frac{7}{4}\zeta \left( 6 \right) + \frac{5}{8}\zeta \left( 2 \right)\zeta \left( 3 \right)\ln 2 - 2\zeta \left( 2 \right){\rm Li}{_4}\left( {\frac{1}{2}} \right)  \\ &\quad\ \ \ \ \ \ \ \  \ \ \ \ \ \ \ \ \ \ \ \ \ \ \ \ \ \ \ + \frac{5}{4}\zeta \left( 4 \right){\ln ^2}2- \frac{1}{{12}}\zeta \left( 2 \right){\ln ^4}2,\tag{3.4}\\
&S(2,3;\bar 1)
 = - \frac{{161}}{{64}}\zeta \left( 6 \right) + \frac{{31}}{{16}}\zeta \left( 5 \right)\ln 2 + \frac{9}{{32}}{\zeta ^2}\left( 3 \right) + \frac{3}{8}\zeta \left( 2 \right)\zeta \left( 3 \right)\ln 2 + 2\zeta \left( 2 \right){\rm Li}{_4}\left( {\frac{1}{2}} \right) \\
&\quad \ \ \ \ \ \ \ \ \ \ \ \ \  - \frac{5}{4}\zeta \left( 4 \right){\ln ^2}2 + \frac{1}{{12}}\zeta \left( 2 \right){\ln ^4}2 + S(2;\bar 4)  - S(\bar 3;3).\tag{3.5}
\end{align*}
In [21], we gave the following formula
\begin{align*}
 S(\bar 1,2;3)&=\frac{{29}}{8}\zeta \left( 2 \right)\zeta \left( 3 \right)\ln 2 - \frac{{93}}{{32}}\zeta \left( 5 \right)\ln 2 - \frac{{1855}}{{128}}\zeta \left( 6 \right) + \frac{{17}}{{16}}{\zeta ^2}\left( 3 \right) \\
& \quad - S(\bar 1;\bar 5) + S(\bar 2;4)+ 4S(2;\bar 4) + 8S(1;\bar 5).\tag{3.6}
\end{align*}
Substituting (3.6) into (3.4) respectively, we arrive at the conclusion that
\begin{align*}
S(\bar 1,3;2) &=\frac{{2079}}{{128}}\zeta \left( 6 \right) + \frac{{93}}{{32}}\zeta \left( 5 \right)\ln 2 - \frac{5}{{16}}{\zeta ^2}\left( 3 \right) - 3\zeta \left( 2 \right)\zeta \left( 3 \right)\ln 2 - 2\zeta \left( 2 \right){\rm Li}{_4}\left( {\frac{1}{2}} \right) \\
& \quad + \frac{5}{4}\zeta \left( 4 \right){\ln ^2}2 - \frac{1}{{12}}\zeta \left( 2 \right){\ln ^4}2 + S(\bar 1;\bar 5)  - S(\bar 2;4) - 4S(2;\bar 4) - 8S(1;\bar 5).\tag{3.7}
\end{align*}
Similarly, taking $\left( {x,y,z} \right) = \left( {1,1,1} \right),\left( {{l_1},{l_2},m} \right) = \left( {2,p + 2m + 1,p} \right)$ in (2.19), we get
\begin{align*}
&S\left( {2,p + 2m + 1;p} \right) + S\left( {2,p;p + 2m + 1} \right) + S\left( {p,p + 2m + 1;2} \right) \\
& =S\left( {2;2p + 2m + 1} \right) + S\left( {p;p + 2m + 3} \right) + S\left( {p + 2m + 1;p + 2} \right)\\
&\quad + \zeta \left( 2 \right)\zeta \left( p \right)\zeta \left( {p + 2m + 1} \right) - \zeta \left( {2p + 2m + 3} \right) . \tag{3.8}
\end{align*}
From (3.1) and (3.8), we know that the quadratic sums
\[S(p,p+2m+1;2)=\sum\limits_{n = 1}^\infty  {\frac{{{\zeta _n}\left( p \right){\zeta _n}\left( {p + 2m + 1} \right)}}{{{n^2}}}} \]
can be expressed in terms of zeta values, linear sums and
\[S(1,p;p+2m+2)=\sum\limits_{n = 1}^\infty  {\frac{{{H_n}{\zeta _n}\left( p \right)}}{{{n^{p + 2m + 2}}}}} ,\ S(1,p+2m+1;p+1)=\sum\limits_{n = 1}^\infty  {\frac{{{H_n}{\zeta _n}\left( {p + 2m + 1} \right)}}{{{n^{p + 1}}}}}. \]
In the last of this section, we give some examples. First, in \cite{X2016}, we showed that all quadratic Euler sums of the form
\[S({1,m;p})=\sum\limits_{n = 1}^\infty  {\frac{{{H_n}{\zeta _n}\left( m \right)}}{{{n^p}}}},\ \ \left( {m + p+1 \le 9} \right)\]
are reducible to $\mathbb{Q}$-linear combinations of single zeta monomial with the addition of linear sums \{S(2;6)\} for weight 8 and give explicit formulas.
From (3.1), (3.2) and (3.8), we can give the following examples
\begin{align*}
&   {2S(2,3;2)+S(2,2;3) = \frac{{13}}{2}\zeta \left( 3 \right)} \zeta \left( 4 \right) - 3\zeta \left( 7 \right),\\
&  S(2,3;2)+S(2,2;3)=  - \frac{{179}}{{16}}\zeta \left( 7 \right) + 8\zeta \left( 3 \right)\zeta \left( 4 \right)+ \frac{5}{2}\zeta \left( 2 \right)\zeta \left( 5 \right),\\
&S(2,2;4)+2S(2,4;2)  = 3\zeta \left( 8 \right) + 2{S({2;6})},\\
&S(2,2;4)-S(2,4;2)  = \frac{{1317}}{{36}}\zeta \left( 8 \right) - 60\zeta \left( 3 \right)\zeta \left( 5 \right) + 9\zeta \left( 2 \right){\zeta ^2}\left( 3 \right) + \frac{{31}}{2}{S({2;6})},\\
&2S(2,5;2)+S(2,2;5) = \frac{{55}}{2}\zeta \left( 9 \right) - 21\zeta \left( 2 \right)\zeta \left( 7 \right) + 4\zeta \left( 3 \right)\zeta \left( 6 \right) + \frac{{13}}{2}\zeta \left( 4 \right)\zeta \left( 5 \right),\\
&S(2,5;2)+S(2,2;5) =  - \frac{{79}}{{72}}\zeta \left( 9 \right) - 7\zeta \left( 2 \right)\zeta \left( 7 \right) + \frac{4}{3}\zeta \left( 3 \right)\zeta \left( 6 \right) + \frac{{23}}{2}\zeta \left( 4 \right)\zeta \left( 5 \right) + \frac{2}{3}{\zeta ^3}\left( 3 \right),\\
&S(2,4;3)+S(2,3;4)  =  - 35\zeta \left( 9 \right) + 14\zeta \left( 2 \right)\zeta \left( 7 \right) + \frac{{107}}{{12}}\zeta \left( 3 \right)\zeta \left( 6 \right) + \frac{7}{2}\zeta \left( 4 \right)\zeta \left( 5 \right) - \frac{1}{3}{\zeta ^3}\left( 3 \right),\\
&S(2,4;3)+S(2,3;4)+S(3,4;2) = - \frac{{77}}{2}\zeta \left( 9 \right) + 21\zeta \left( 2 \right)\zeta \left( 7 \right) + \frac{{15}}{4}\zeta \left( 3 \right)\zeta \left( 6 \right) + 3\zeta \left( 4 \right)\zeta \left( 5 \right).
\end{align*}
Therefore, combining related equations, we obtain the following identities.
\begin{exa} Some results on quadratic Euler sums.
\begin{align*}
&S(2,2;3)=\sum\limits_{n = 1}^\infty  {\frac{{\zeta _n^2\left( 2 \right)}}{{{n^3}}}}  =  - \frac{{155}}{8}\zeta \left( 7 \right) + \frac{{19}}{2}\zeta \left( 3 \right)\zeta \left( 4 \right) + 5\zeta \left( 2 \right)\zeta \left( 5 \right),\\
&S(2,3;2)=\sum\limits_{n = 1}^\infty  {\frac{{{\zeta _n}\left( 2 \right){\zeta _n}\left( 3 \right)}}{{{n^2}}}}  = \frac{{131}}{{16}}\zeta \left( 7 \right) - \frac{3}{2}\zeta \left( 3 \right)\zeta \left( 4 \right) - \frac{5}{2}\zeta \left( 2 \right)\zeta \left( 5 \right),\\
&S(2,2;4)=\sum\limits_{n = 1}^\infty  {\frac{{\zeta _n^2\left( 2 \right)}}{{{n^4}}}}  = 11{S({2;6})} + \frac{{457}}{{18}}\zeta \left( 8 \right) + 6\zeta \left( 2 \right){\zeta ^2}\left( 3 \right) - 40\zeta \left( 3 \right)\zeta \left( 5 \right), \\
&S(2,4;2)=\sum\limits_{n = 1}^\infty  {\frac{{{\zeta _n}\left( 2 \right){\zeta _n}\left( 4 \right)}}{{{n^2}}}}  =  - \frac{9}{2}{S({2;6})} - \frac{{403}}{{36}}\zeta \left( 8 \right) - 3\zeta \left( 2 \right){\zeta ^2}\left( 3 \right) + 20\zeta \left( 3 \right)\zeta \left( 5 \right),\\
&S(2,2;5)=\sum\limits_{n = 1}^\infty  {\frac{{\zeta _n^2\left( 2 \right)}}{{{n^5}}}}  =  - \frac{{1069}}{{36}}\zeta \left( 9 \right) + \frac{4}{3}{\zeta ^3}\left( 3 \right) + 7\zeta \left( 2 \right)\zeta \left( 7 \right) - \frac{4}{3}\zeta \left( 3 \right)\zeta \left( 6 \right) + \frac{{33}}{2}\zeta \left( 4 \right)\zeta \left( 5 \right),\\
&S(2,5;2)=\sum\limits_{n = 1}^\infty  {\frac{{{\zeta _n}\left( 2 \right){\zeta _n}\left( 5 \right)}}{{{n^2}}}}  = \frac{{2059}}{{72}}\zeta \left( 9 \right) - \frac{2}{3}{\zeta ^3}\left( 3 \right) - 14\zeta \left( 2 \right)\zeta \left( 7 \right) + \frac{8}{3}\zeta \left( 3 \right)\zeta \left( 6 \right) - 5\zeta \left( 4 \right)\zeta \left( 5 \right),\\
&S(3,4;2)=\sum\limits_{n = 1}^\infty  {\frac{{{\zeta _n}\left( 3 \right){\zeta _n}\left( 4 \right)}}{{{n^2}}}}  =  - \frac{7}{2}\zeta \left( 9 \right) + 7\zeta \left( 2 \right)\zeta \left( 7 \right) - \frac{{31}}{6}\zeta \left( 3 \right)\zeta \left( 6 \right) - \frac{1}{2}\zeta \left( 4 \right)\zeta \left( 5 \right) + \frac{1}{3}{\zeta ^3}\left( 3 \right).
\end{align*}
\end{exa}
Moreover, we use Mathematica tool to check numerically each of the specific identities listed. The numerical values of
nonlinear Euler sums of weights \{7,8,9,10\}, to 30 decimal digits, see Table 1.
\newcommand{\tabincell}[2]{\begin{tabular}{@{}#1@{}}#2\end{tabular}}
\begin{table}[htbp]
 \centering
\caption{\label{tab:test}Numerical approximation}
 \begin{tabular}{|c|c|c|}
  \hline
 \tabincell{c}{Euler sum} &\tabincell{c}{ Numerical values of closed form\\ (30 decimal digits)}
  & \tabincell{c}{Numerical approximation of\\ Euler sum (30 decimal digits)}\\
  \hline
  $S({2,2;3})$ &1.35125578526281388688070479101&1.35125578526281388688070478635 \\
 \hline
  $S({2,3;2})$ &2.04014406352629668230178759593&2.04014406352629668230178759172 \\
\hline
$S({1,2;5})$ &1.07388087034296588059339568891&1.07388087034296588059339568663 \\
\hline
$S({1,3;4})$& 1.15201859049597540982393939989& 1.15201859049597540982393939372\\
\hline
$S({1,4;3})$& 1.37755320390542981268777869872 &1.37755320390542981268777869712\\
\hline
$S({1,5;2})$&2.45339834780017683307966649793 & 2.45339834780017683307966649461\\
\hline
$S({2,2;4})$& 1.13642391274089928376327915373 &1.13642391274089928376327915559\\
\hline
$S({2,4;2})$&1.95980117454124719492773304920 &1.95980117454124719492773304287\\
\hline
$S({2,2;5})$&1.05972458873705638208576920975&1.05972458873705638208576920818\\
\hline
$S({2,5;2})$&1.92499254625584068819896689186&1.92499254625584068819896688762\\
\hline
$S({3,4;2})$&1.80313006078587093607835773253&1.80313006078587093607835772809\\
\hline
$S({1,2;7})$&1.01603499621822946463309621255&1.01603499621822946463309621221\\
\hline
$S({1,3;6})$&1.03017876630576928913732006061&1.03017876630576928913732005893\\
\hline
$S({1,4;5})$&1.06164990978502285301181351196&1.06164990978502285301181351270\\
\hline
$S({1,5;4})$&1.13783419529420067466663388885&1.13783419529420067466663388537\\
\hline
$S({1,6;3})$&1.35867450783449320806721637607&1.35867450783449320806721637012\\
\hline
$S({1,7;2})$&2.41561649536052525591387317796&2.41561649536052525591387317514\\
\hline
 \end{tabular}
\end{table}
\\
In fact, by using the methods of this paper, it is possible to establish other identities of Euler sums. For example, taking $m=1$ in (2.11), we deduce that
\begin{align*}
\int\limits_0^1 {{x^{n - 1}}\ln x{\rm{L}}{{\rm{i}}_p}\left( x \right)} dx =& \sum\limits_{i = 1}^{p - 1} {\sum\limits_{j = 1}^{p - i} {{{\left( { - 1} \right)}^{i + j - 1}}\frac{{\zeta \left( {p + 2 - i - j} \right)}}{{{n^{i + j}}}}} } \\
&\quad+ {\left( { - 1} \right)^p}p\frac{{{H_n}}}{{{n^{p + 1}}}} + {\left( { - 1} \right)^p}\frac{{{\zeta _n}\left( 2 \right) - \zeta \left( 2 \right)}}{{{n^p}}}. \tag{3.9}
\end{align*}
Multiplying (2.7) and (2.16) by $\frac{\ln{x}}{x}$, and integrating over the interval (0,1), we arrive at the conclusion that
\begin{align*}
&{\left( { - 1} \right)^p}\left[ {S\left( {3,p + 2m + 1;p} \right) + S\left( {3,p;p + 2m + 1} \right)} \right]\\
 = &{\left( { - 1} \right)^p}\zeta \left( 3 \right)\left[ {S\left( {p + 2m + 1;p} \right) + S\left( {p;p + 2m + 1} \right)} \right]\\
 & + {\left( { - 1} \right)^{p + 1}}\frac{{p\left( {p + 1} \right)}}{2}S\left( {1,p + 2m + 1;p + 2} \right)\\
 & + {\left( { - 1} \right)^{p + 1}}\frac{{\left( {p + 2m + 1} \right)\left( {p + 2m} \right)}}{2}S\left( {1,p;p + 2m + 3} \right)\\
 & + {\left( { - 1} \right)^{p + 1}}p\left[ {S\left( {2,p + 2m + 1;p + 1} \right) - \zeta \left( 2 \right)S\left( {p + 2m + 1;p + 1} \right)} \right]\\
 & + {\left( { - 1} \right)^{p + 1}}\left( {p + 2m + 1} \right)\left[ {S\left( {2,p;p + 2m + 2} \right) - \zeta \left( 2 \right)S\left( {p;p + 2m + 2} \right)} \right]\\
 &- \sum\limits_{l = 1}^{p - 1} {\sum\limits_{i = 1}^{p - l} {\sum\limits_{j = 1}^{p + 1 - i - l} {{{\left( { - 1} \right)}^{i + j + l}}\zeta \left( {p + 3 - i - j - l} \right)S\left( {p + 2m + 1;i + j + l} \right)} } } \\
 & + \sum\limits_{l = 1}^{p + 2m} {\sum\limits_{i = 1}^{p + 2m + 1 - l} {\sum\limits_{j = 1}^{p + 2m + 2 - i - l} {{{\left( { - 1} \right)}^{i + j + l}}\zeta \left( {p + 2m + 4 - i - j - l} \right)S\left( {p;i + j + l} \right)} } } ,\tag{3.10}
\end{align*}
and
\begin{align*}
 &{\left( { - 1} \right)^p}\left[ {S\left( {3,p + 2m;p} \right) - S\left( {3,p;p + 2m} \right)} \right] \\
  =& {\left( { - 1} \right)^p}\zeta \left( 3 \right)\left[ {S\left( {p + 2m;p} \right) - S\left( {p;p + 2m} \right)} \right] \\
 & + {\left( { - 1} \right)^{p + 1}}\frac{{p\left( {p + 1} \right)}}{2}S\left( {1,p + 2m;p + 2} \right) \\
  &+ {\left( { - 1} \right)^p}\frac{{\left( {p + 2m} \right)\left( {p + 2m + 1} \right)}}{2}S\left( {1,p;p + 2m + 2} \right) \\
  &+ {\left( { - 1} \right)^{p + 1}}p\left[ {S\left( {2,p + 2m;p + 1} \right) - \zeta \left( 2 \right)S\left( {p + 2m;p + 1} \right)} \right] \\
  &+ {\left( { - 1} \right)^p}\left( {p + 2m} \right)\left[ {S\left( {2,p;p + 2m + 1} \right) - \zeta \left( 2 \right)S\left( {p;p + 2m + 1} \right)} \right] \\
  &- \sum\limits_{l = 1}^{p - 1} {\sum\limits_{i = 1}^{p - l} {\sum\limits_{j = 1}^{p + 1 - i - l} {{{\left( { - 1} \right)}^{i + j + l}}\zeta \left( {p + 3 - i - j - l} \right)S\left( {p + 2m;i + j + l} \right)} } }  \\
  &+ \sum\limits_{l = 1}^{p + 2m - 1} {\sum\limits_{i = 1}^{p + 2m - l} {\sum\limits_{j = 1}^{p + 2m + 1 - i - l} {{{\left( { - 1} \right)}^{i + j + l}}\zeta \left( {p + 2m + 3 - i - j - l} \right)S\left( {p;i + j + l} \right)} } }.  \tag{3.11}
\end{align*}
Putting $m=1,p=2$ in (3.11), we have
\[\sum\limits_{n = 1}^\infty  {\left\{ {\frac{{{\zeta _n}\left( 3 \right){\zeta _n}\left( 4 \right)}}{{{n^2}}} - \frac{{{\zeta _n}\left( 3 \right){\zeta _n}\left( 2 \right)}}{{{n^4}}}} \right\}}  =  - \frac{{1063}}{{36}}\zeta \left( 9 \right) + 14\zeta \left( 2 \right)\zeta \left( 7 \right) - \frac{{37}}{6}\zeta \left( 3 \right)\zeta \left( 6 \right) + 13\zeta \left( 4 \right)\zeta \left( 5 \right).\]
Combining related equations, we can obtain
\begin{align*}
 &\sum\limits_{n = 1}^\infty  {\frac{{{\zeta _n}\left( 2 \right){\zeta _n}\left( 3 \right)}}{{{n^4}}}}  = \frac{{937}}{{36}}\zeta \left( 9 \right) + \frac{1}{3}{\zeta ^3}\left( 3 \right) - 7\zeta \left( 2 \right)\zeta \left( 7 \right) + \zeta \left( 3 \right)\zeta \left( 6 \right) - \frac{{27}}{2}\zeta \left( 4 \right)\zeta \left( 5 \right), \\
 &\sum\limits_{n = 1}^\infty  {\frac{{{\zeta _n}\left( 2 \right){\zeta _n}\left( 4 \right)}}{{{n^3}}}}  =  - \frac{{2197}}{{36}}\zeta \left( 9 \right) - \frac{2}{3}{\zeta ^3}\left( 3 \right) + 21\zeta \left( 2 \right)\zeta \left( 7 \right) + \frac{{95}}{{12}}\zeta \left( 3 \right)\zeta \left( 6 \right) + 17\zeta \left( 4 \right)\zeta \left( 5 \right).
\end{align*}
{\bf Acknowledgments.} We thank the anonymous referee for suggestions which led to improvements in the exposition.


{\small
\begin{thebibliography}{99}

\bibitem{A2000} George E. Andrews, Richard Askey, Ranjan Roy. {\sl Special Functions}.
Cambridge University Press., 2000: 481-532.
\bibitem{BBG1994}
David H. Bailey, Jonathan M. Borwein and Roland Girgensohn. {\sl Experimental evaluation of Euler sums}.
Experimental Mathematics., 1994, {\bf 3}(1): 17-30.

\bibitem{BBC2014}
David H. Bailey, Jonathan M. Borwein, Richard E. Crandall. {\sl Computation and theory of extended Mordell-Tornheim-Witten sums}. Math. Comp., 2014, {\bf 83}(288): 1795-1821.
\bibitem{B1985} B. C. Berndt. {\sl Ramanujan¡¯s Notebooks, Part I}. Springer-Verlag, New York., 1985.
\bibitem{B1989} B. C. Berndt. {\sl Ramanujan¡¯s Notebooks, Part II}. Springer-Verlag, New York., 1989.
\bibitem{BBG1995}
David Borwein, Jonathan M. Borwein and Roland Girgensohn. {\sl Explicit evaluation of
Euler sums}. Proc. Edinburgh Math., 1995, {\bf 38}: 277-294.

\bibitem{BBGP1996}
J.Borwein, P.Borwein, R.Girgensohn, S.Parnes. {\sl Making sense of experimental mathematics}. Mathematical Intelligencer., 1996, {\bf18}(4): 12-18.

\bibitem{BBBL2001}
Jonathan M. Borwein, David M. Bradley, David J. Broadhurst, Petr. Lison¨§k.
{\sl Special values of multiple polylogarithms. Trans}. Amer. Math. Soc., 2001, {\bf 353}(3): 907-941.
\bibitem{BZB2008}
J. M. Borwein, I. J. Zucker, J. Boersma. {\sl The evaluation of character Euler double sums}. Ramanujan J., 2008, {\bf 15} (3): 377-405.

\bibitem{BG1996}
J.M. Borwein, R. Girgensohn, {\sl Evaluation of triple Euler sums}, Electron. J. Combin., 1996: 2-7.

\bibitem{EW2012}
 Minking Eie, Chuan-Sheng Wei. {\sl Evaluations of some quadruple Euler sums of even weight}. Functions et Approximatio., 2012, {\bf 46}(1): 63-67.

\bibitem{FS1998}
Philippe Flajolet and Bruno Salvy. {\sl Euler sums and contour integral representations}. Experimental Mathematics., 1998, {\bf 7}(1): 15--35.

\bibitem{F2005}
 Pedro freitas. {\sl Integrals of polylogarithmic functions, recurrence relations, and associated Euler sums}. Mathematics of Computation., 2005, {\bf 74}(251): 1425-1440.


\bibitem{L1974}
Comtet L. Advanced combinatorics, Boston: D Reidel Publishing Company., 1974.

\bibitem{M2014}
I.Mez$\ddot{o}$. {\sl Nonlinear Euler sums}. Pacific J. Math., 2014, {\bf 272}: 201-226.

\bibitem {S2015}
A. Sofo. {\sl Quadratic alternating harmonic number sums}. J. Number Theory., 2015, {\bf 154}: 144-159.
\bibitem {M1994} C. Markett. {\sl Triple sums and the Riemann zeta function}. J. Number Theory., 1994, {\bf 48}(2): 113-132.
\bibitem{S2007}
Ping Sun. {\sl The 6-order Sums of Riemann Zeta Function}. Acta mathematica sinica,chinese Series ., 2007, {\bf 50}(2): \bibitem{Xu2016}
Ce Xu, Jinfa Cheng. {\sl Some Results On Euler Sums}. Functions et Approximatio., 2016, {\bf 54}(1): 25-37.
\bibitem{X2016}
Ce Xu, Yuhuan Yan, Zhijuan Shi. {\sl Euler sums and integrals of polylogarithm functions}. J. Number Theory., 2016, {\bf 165}: 84-108.
\bibitem{Xyz2016}
Ce Xu, Yingyue Yang, Jianwen Zhang. {\sl Explicit evaluation of quadratic Euler sums}. Int. J. Number Theory., 2017, {\bf 13}(2): 655-672.
\bibitem{DZ2012}
D. Zagier. {\sl Evaluation of the multiple zeta values $\zeta(2,...,2,3,2,...,2)$}. Annals of Mathematics., 2012, {\bf 2}(2): 977-1000.
\end{thebibliography}}
 \end{document}